\documentclass[a4paper,11pt]{amsart}

\usepackage{amsmath}
\usepackage{amsthm}
\usepackage{amsfonts}
\usepackage{amssymb}
\usepackage[all]{xy}

\newcommand{\myAA}{\mathcal{A}}
\newcommand{\myFF}{\mathcal{F}}

\newcommand{\ZZ}{\mathbb{Z}}

\newcommand{\QQ}{\mathbb{Q}}
\newcommand{\BB}{\mathbf{B}}

\newtheorem{theorem}{Theorem}[section]
\newtheorem{proposition}{Proposition}[theorem]
\newtheorem{corollary}{Corollary}[theorem]
\newtheorem{definition}{Definition}[theorem]
\newtheorem{lemma}{Lemma}[theorem]

\begin{document}

\title[Quiver Grassmannians and cluster algebras]{Geometry of quiver Grassmannians of Kronecker type and canonical basis of cluster algebras}
\author{Giovanni Cerulli Irelli\and Francesco Esposito}
\footnote{Research supported by grant CPDA071244/07 of Padova University}

\begin{abstract}
We study quiver Grassmannians associated with indecomposable representations of the Kronecker quiver. We find a cellular decomposition of them and we compute their Betti numbers. As an application, we find a geometric realization of the ``canonical basis' '  of cluster algebras of type $A_1^{(1)}$ found by Sherman and Zelevinsky in \cite{Sherman} and of type $A_2^{(1)}$ found in \cite{CanonicalBasis}. 
\end{abstract}
\maketitle

\section{Introduction}\label{Sec:Intro}
Cluster algebras are commutative $\ZZ$--subalgebras of the field of rational functions in a finite number of indeterminates which have been introduced and studied by S.~Fomin and A.~Zelevinsky in a series of papers \cite{FZI}, \cite{FZII}, \cite{BFZIII} and \cite{FZIV}. To every quiver $Q$ without loops and $2$--cycles it is associated a coefficient--free cluster algebra $\myAA_Q$. In \cite{CC}, \cite{CK1}, \cite{CK2} and \cite{DWZII} the authors describe the cluster variables of $\myAA_Q$ via a map, called the Caldero--Chapoton map,   between the representations of $Q$  and the field of rational functions in $n$ variables  (we address the reader to the survey \cite{keller-2008}). Such map is defined in terms of Euler--Poincar\'e characteristic of some complex projective varieties attached to every representation $M$ of $Q$ and called quiver Grassmannians. By definition the quiver Grassmannian $Gr_{\mathbf{e}}(M)$ consists of all sub--representations of $M$ of dimension vector $\mathbf{e}$. These varieties are considered in several places, e.g. \cite{CK1}, \cite{CK2}, \cite{CR}, \cite[section~12.3]{hernandez}, \cite{Nakajima}, and in this paper we try to add some more geometric information about them, at least in the case of the Kronecker quiver. In \cite{CZ} authors compute the Euler--Poincar\'e characteristic of quiver Grassmannians associated with the Kronecker quiver and they conjecture the existence of a cellular decomposition which we find here. In \cite{String} a torus action on some quiver Grassmannians has been found and this allows to produce a cellular decomposition of them in the case they are smooth.

In this paper we study quiver Grassmannians associated with the Kronecker quiver. A representation of the Kronecker quiver (in the sequel we say a $Q$--representation) is a quadruple $M=(M_1, M_2,m_a,m_b)$ where $M_1$ and $M_2$ are finite dimensional complex vector spaces and $m_a,m_b:M_1\rightarrow M_2$ are two linear maps between them. Given two non--negative integers $e_1$ and $e_2$, the variety $Gr_{(e_1,e_2)}(M)$ is defined as the set
$$
\{(N_1,N_2)\in Gr_{e_1}(M_1)\times Gr_{e_2}(M_2):\; m_a(N_1)\subset N_2,\;m_b(N_1)\subset N_2\}
$$
where $Gr_e(V)$ denotes the Grassmannians of $e$--dimensional vector spaces in a vector space $V$. This is a projective variety which is in general not smooth. When $M_1=M_2=\mathbf{C}^n$ and $m_a=Id$ is the identity and $m_b=J_n(0)$ is an indecomposable nilpotent Jordan block, the representation $M$ is regular indecomposable and we denote it by $R_n$. The corresponding quiver Grassmannians $Gr_{(e_1,e_2)}(R_n)$ are the main subjects of this paper. We briefly write $X$ for one of these quiver Grassmannians. We find that a one--dimensional torus $T$ acts on $X$ (section~\ref{Sec:ActionOnX}). We provide a stratification of $X$ (see section~\ref{Sec:Stratification})
\begin{equation}\label{Eq:StratificationIntro}
X=X_0\supseteq X_1\supseteq\cdots \supseteq X_s
\end{equation}
for $s=\textrm{min}(e_1,n-e_2)$ into closed subvarieties $X_k \simeq Gr_{(e_1-k,e_2-k)}(R_{n-2k})$. Moreover $X_{k+1}$ is the singular locus of $X_k$ and the difference $X_k\setminus X_{k+1}$ is  a smooth quasi--projective variety which is not complete. As a consequence of the stratification \eqref{Eq:StratificationIntro} we get that $X$ is smooth if and only if $s=0$, i.e. either $e_1=0$ or $e_2=n$, in which cases the quiver Grassmannian specializes to an usual Grassmannian of vector subpaces.

In section~\ref{Sec:CellDec} we prove that Bia{\l}ynicki-Birula's theorem on cellular decomposition of smooth projective varieties applies to $X_k\setminus X_{k+1}$ and we can hence prove that $X$ decomposes
$$
X=\cup_{L\in X^T}X_L
$$
into attracting sets of $T$--fixed points
$$
X_L:=\{N\in X| \lim_{\lambda\rightarrow 0}t_\lambda N=L\}
$$
and theses sets are affine spaces. In section~\ref{Sec:Cells} we describe the cell $X_L$: if $L$ is indecomposable then $X_L$ is the orbit of $L$ by the action of a  unipotent group; if $L=L'\oplus L''$ with $L'$ of ``lower weight than'' $L''$ then the following formula holds:
$$
\textrm{dim}X_{L'\oplus L''}=\textrm{dim} X_{L'}+\textrm{dim}X_{L'}-\langle\mathbf{dim}L',\mathbf{dim}L''\rangle
$$
where $\langle\cdot,\cdot\rangle$ denotes the Euler form of the Kronecker quiver. As a consequence of this formula we are able to compute the Poincar\'e polynomials of every indecomposable representation
of the Kronecker quiver: let $P_n$ (resp. $I_n$) be the indecomposable preprojective of dimension $(n,n+1)$ (resp. $(n+1,n)$). We have (see section~\ref{Sec:Betti})
\begin{equation}\label{Eq:IntroPoincare}
\begin{array}{ccc}
P_{Gr_{\mathbf{e}}(R_n)}(t)&=&P_{Gr_{(e_2-e_1)}(e_2)}(t)P_{Gr_{(e_2-e_1)}(n-e_1)}(t)\\
P_{Gr_\mathbf{e}(P_n)}(t)&=&P_{Gr_{e_1}(e_2-1)}(t)P_{Gr_{(e_2-e_1)}(n+1-e_1)}(t)
\\
P_{Gr_\mathbf{e}(I_n)}(t)&=&P_{Gr_{e_1}(e_2+1)}(t)P_{Gr_{(e_2-e_1)}(n-e_1)}(t)
\end{array}
\end{equation}
where $Gr_t(s)$ denotes the Grassmannian of $t$--dimensional vector subspaces of an $s$--dimensional vector space. 

We remark that formulas \eqref{Eq:IntroPoincare} yield another proof of results of Szanto \cite{szanto-2009}. Indeed 
$Gr_\mathbf{e}(P_n)$ and $Gr_\mathbf{e}(I_n)$ are smooth over any field and hence this follows from \cite{BB}. For $X=Gr_\mathbf{e}(R_n)$, both  $X$ and all the $X_k$'s may be defined over $\ZZ$ and the results of theorem~\ref{Thm:Stratification} continues to hold over $\ZZ$. In particular $X_k \setminus X_{k+1}$ is smooth over $\ZZ$ and by  
standard theorems in $l$-adic cohomology one has that the cohomology  
of the quiver grassmannian over the field $F_q$ is the same as in  
characteristic zero and the frobenius acts by a suitable power of $q$.  
Applying the Lefschetz fixed point formula, our formula for the Poincar\'e polynomials provide another proof of  
the results in \cite{szanto-2009}. It would be interesting to know if the cellular decomposition continues to hold in positive characteristic.

In section~\ref{Sec:ClusterAlgebras} we consider the (coefficient--free) cluster algebra $\myAA$ of type $A_1^{(1)}$. The cluster variables of $\myAA$ are the rational functions $x_m$, $m\in\ZZ$, of the field $\myFF=\QQ(x_1,x_2)$ recursively generated by the following relation:
$$
x_kx_{k+2}=x_{k+1}^2+1.
$$
%It is not hard to see that every couple $\{x_k,x_{k+1}\}$ of consecutive cluster variables is a free--generating system for the field $\myFF$ and hence $\myFF=\QQ(x_k,x_{k+1})$ and every cluster variable can be expressed as a rational function in every such couple (which are called the clusters of $\myAA$). Moreover every cluster variable is a Laurent polynomial in every cluster of $\myAA$ (this is the Laurent phenomenon proved in \cite{FZI}). 
In this case the Caldero--Chapoton map $M\mapsto CC(M)$ associates to a $Q$--representation $M$ of dimension vectors $(d_1,d_2)$ the following Laurent polynomial:
\begin{equation}\label{Eq:CCIntro}
CC(M):=\frac{\sum_{\mathbf{e}}\chi(Gr_\mathbf{e}(M))x_1^{2(d_2-e_2)}x_2^{2e_1}}{x_1^{d_1}x_2^{d_2}}
\end{equation}
In \cite{CK1} it is proved that the map $M\mapsto CC(M)$ restricts to a bijection between the indecomposable rigid $Q$--representations $M$ (i.e. $Ext^1(M,M)=0$) and the cluster variables of $\myAA$ different from $x_1$ and $x_2$. Moreover it has the following property:
$CC(M\oplus N)=CC(M)CC(N)$,
under which cluster monomials not divisible by $x_1$ or $x_2$, i.e. monomials of the form $x_k^ax_{k+1}^b$ for $k\in\ZZ\setminus\{1,2\}$ and $a,b\geq0$, are in bijection with rigid $Q$--representations.

%Following \cite{Sherman} we say that an element of $\myAA$ is \emph{positive} if its Laurent expansion in every cluster of $\myAA$ has positive coefficients. Positive elements form a cone which has been described in 
In \cite{Sherman}  the authors introduce distinguished elements $\{z_n|\; n\geq1\}$  of $\myFF$ recursively defined by:
\begin{equation}\label{Eq:ZnDef}
\begin{array}{rcl}
z_0&=&2\\
z_1&=&x_0x_3-x_1x_2\\
z_{n+1}&=&z_1 z_n-z_{n-1}\;\;\; n\geq1
\end{array}
\end{equation}
and they prove that the set
$$
\BB:=\{\textrm{cluster monomials}\}\cup\{z_n: n\geq1\}
$$
is a $\ZZ$--basis of $\myAA$ such that positive linear combinations of its elements coincide with the set of all \emph{positive} elements of $\myAA$ (i.e. elements which are positive Laurent polynomials in every cluster of $\myAA$). They call this basis the canonical basis of $\myAA$. We give a geometric realization of $\BB$ by using the Caldero--Chapoton map: cluster monomials are image of  rigid representations and quiver Grassmannians associated with rigid quiver representations are smooth \cite{CR}. Having this in mind in section~\ref{Sec:ClusterAlgebras} we prove the following: for every $n\geq1$:
$$
z_n=\frac{\sum_{\mathbf{e}}\chi(Gr_\mathbf{e}(R_n)^{sm})x_1^{2(n-e_2)}x_2^{2e_1}}{x_1^{n}x_2^{n}}
$$
where $Gr_\mathbf{e}(R_n)^{sm}:=X_0\setminus X_1$ denotes the smooth part of $Gr_\mathbf{e}(R_n)$.
%In other words the canonical basis can be obtained by computing the CC--map in  the smooth part of the quiver Grassmannians associated with all the rigid representations and of all regular indecomposable representations. 

A similar construction can be made in a cluster algebra $\myAA_2$ of type $A_2^{(1)}$. These cluster algebras are studied in \cite{CanonicalBasis} and some results are recalled in section~\ref{Sec:A21}: the canonical basis of $\myAA_2$ consists of cluster monomials together with elements $\{u_nw^k,\;u_n z^k:\;k\geq0,\; n\geq1\}$ where $w$ and $z$ are two cluster variables and the $u_n$ are defined similarly to \eqref{Eq:ZnDef} as follows:
\begin{equation}\label{Eq:UnDef}
\begin{array}{rcl}
u_0&=&2\\
u_1&=&zw-2\\
u_{n+1}&=&u_1 u_n-u_{n-1}\;\;\; n\geq1
\end{array}
\end{equation}
We prove that $u_n$ is obtained by evaluating the Caldero--Chapoton map at the \emph{smooth part} of the regular indecomposable representation of a quiver of type $A_2^{(1)}$ of dimension $(1,1,1)$ which lies in an homogeneous tube. 

\section{Geometric structure of quiver Grassmannians}\label{Sec:GeomX}

Let
$$
Q:=\xymatrix{1\ar@<1ex>^a[r]\ar_b[r]&2}
$$
be the Kronecker quiver. As usual, we denote a complex $Q$--representation $M=(M_1,M_2,m_a,m_b)$  as follows:
$$
M=\xymatrix{M_1\ar@<1ex>^{m_a}[r]\ar_{m_b}[r]&M_2}.
$$
A sub--representation $N$ of $M$ consists of vector subspaces $N_1$ and $N_2$ of $M_1$ and $M_2$ respectively such that $m_a(N_1)\subset N_2$ and $m_b(N_1)\subset N_2$. The vector $\mathbf{dim}(M):=(dim M_1,dim M_2)$ is called the dimension vector of $M$. A morphism $g:M\rightarrow M'$ from a $Q$--representation $M$ to a $Q$--representation $M'$ is a couple $(g_1,g_2)$ of linear maps $g_1:M_1\rightarrow M_1'$ and $g_2:M_2\rightarrow M_2'$ such that $m_a'\circ g_1=g_2\circ m_a$ and $m'_{b}\circ g_1=g_2\circ m_b$.
% , i.e. the following diagram commutes:
% \begin{equation}\label{Eq:DiagramMorphism}
% \xymatrix{
% M_1\ar_{g_1}[d]\ar@<1ex>^{m_a}[r]\ar_{m_b}[r]&M_2\ar^{g_2}[d]\\
% M'_1\ar@<1ex>^{m'_a}[r]\ar_{m'_b}[r]&M'_2
% }
% \end{equation}
The set of $Q$--representations form a category which is an abelian Krull--Schmidt category via the natural notions of direct sums, kernel and cokernel (see e.g. \cite{ASS}). The classification of $Q$--representations which are \emph{indecomposable} (i.e. that are not direct sum of two non--trivial sub--representations) goes back to Kronecker \cite{Kronecker}. The following is their complete list: there are the \emph{indecomposable preprojectives} (for $n\geq0$):
$$
P_n=\xymatrix{k^n\ar@<1ex>^{\varphi_1}[r]\ar_{\varphi_2}[r]&k^{n+1}}
$$
where $k=\mathbf{C}$ denotes the field of complex numbers, $\varphi_1$, $\varphi_2$ are the two immersions respectively in the vector subspace spanned by the first and the last basis vectors. There are the \emph{indecomposable regulars}:
$$
R_n(\lambda)=\xymatrix{k^n\ar@<1ex>^{Id}[r]\ar_{J_n(\lambda)}[r]&k^{n}},\;\;\;R_n(\infty)=\xymatrix{k^n\ar@<1ex>^{J_n(0)}[r]\ar_{Id}[r]&k^{n}}
$$
where $\lambda\in k$ and $J_n(\lambda)$ denotes the $n-th$ indecomposable Jordan block of eigenvalue $\lambda$ and $Id$ denotes the identity matrix. And finally there are the \emph{indecomposable preinjectives} (for $n\geq0$):
$$
I_n=\xymatrix{k^{n+1}\ar@<1ex>^{\varphi_2^t}[r]\ar_{\varphi^t_1}[r]&k^{n}}
$$
where $\varphi_1^t$ and $\varphi_2^t$ are the transpose of the matrices $\varphi_1$ and $\varphi_2$ defined above.  For an indecomposable $Q$--representation $M$ we denote by $\mathcal{B}_i=\{v_k^{(i)}\}$, $i=1,2$, the standard basis of $M_i$ with respect to which $M$ has the previous presentation. All the other $Q$--representations are direct sums of these ones. Direct sums of indecomposable preprojectives (resp. regulars, preinjectives) are called preprojective (resp. regular, preinjective) $Q$--representations.

Given non--negative integers $e_1$ and $e_2$ and a $Q$--representation $M$ we consider the variety:
$$
Gr_\mathbf{e}(M):=\{N\leq M|\; \mathbf{dim}(N)=(e_1,e_2)\}
$$
called the quiver Grassmannian of dimension $\mathbf{e}=(e_1,e_2)$ of $M$. This is closed inside the product $Gr_{e_1}(M_1)\times Gr_{e_2}(M_2)$ of usual Grassmannians of vector subspaces and it is hence a complex projective variety. In \cite{CR} it is shown that the tangent space $T_N(Gr_{\mathbf{e}}(M))$ at a point $N$ of $Gr_{\mathbf{e}}(M)$ equals:
\begin{equation}\label{Eq:TN}
T_N(Gr_{\mathbf{e}}(M))=Hom(N,M/N).
\end{equation}
Moreover the following inequalities hold for $Z:=Gr_\mathbf{e}(M)$:
\begin{equation}
\langle\mathbf{e},\mathbf{d}-\mathbf{e}\rangle\leq \textrm{dim}Z\leq \textrm{dim}T_N(Z)\leq \langle\mathbf{e},\mathbf{d}-\mathbf{e}\rangle+\textrm{dim} Ext^1(M,M)
\end{equation}
where $\langle(a,b)^t,(c,d)^t\rangle:=ac+bd-2ad$ is the Euler form of $Q$ and $\mathbf{d}:=\mathbf{dim}M$ so that $\mathbf{dim}(M/N)=\mathbf{d}-\mathbf{e}$. In particular if $M$ is \emph{rigid}, i.e. $Ext^1(M,M)=0$, then all the quiver Grassmannians associated with it are smooth (see also \cite[proposition~3.5]{DWZII} for another such result) and they have dimension $\textrm{dim}Gr_\mathbf{e}(M)=\langle\mathbf{e},\mathbf{d}-\mathbf{e}\rangle$. It is known that the rigid $Q$--representations are the following:
$$
P_n^{\oplus a}\oplus P_{n+1}^{\oplus b},\;\;\;I_n^{\oplus a}\oplus I_{n+1}^{\oplus b}
$$
for all $n\geq0$ and $a,b\geq0$.

\subsection{Action of a group on quiver Grassmannians}\label{Sec:GroupAction}
Let $M$ be  a $Q$--representation of dimension vector $\mathbf{d}=(d_1,d_2)$ we consider the group
$$
\mathcal{N}(M):=
\{(A,B)\in \prod_{i=1}^{2}GL_{d_i}(M_i)|\exists\; \lambda\in\mathbf{C}^\ast: m_aA=Bm_a,\;m_bA=\lambda Bm_b\}.
$$
Note that the automorphism group of $M$ is a closed subgroup of $\mathcal{N}(M)$ (for $\lambda=1$). The group $\mathcal{N}(M)$ acts on $Gr_\mathbf{e}(M)$ as follows $(A,B)\cdot(N_1,N_2):=(AN_1,BN_2)$.
\begin{lemma}
\begin{enumerate}
\item For every $n\geq0$, $\mathcal{N}(P_n)\simeq \mathbf{C}^\ast\times \mathbf{C}^\ast$;
\item For every $n\geq0$, $\mathcal{N}(I_n)\simeq \mathbf{C}^\ast\times \mathbf{C}^\ast$.
\end{enumerate}
\end{lemma}
\begin{proof}
It follows easily from the definition that $\mathcal{N}(P_n)$ consists of diagonal matrices $(A,B)$ of the form:
$$
A=\textrm{diag}(a,a\lambda,\cdots, a\lambda^{n-1})
$$ 
$$
B=\textrm{diag}(a,a\lambda,\cdots, a\lambda^{n})
$$ 
for $a,\lambda\in\mathbf{C}^\ast$. Similarly for $\mathcal{N}(I_n)$.
\end{proof}

\begin{proposition}
Let $M$ be a rigid $Q$--representation. Then every quiver Grassmannian $Gr_\mathbf{e}(M)$ associated with $M$ has a cellular decomposition.
\end{proposition}
\begin{proof}
The variety $Gr_\mathbf{e}(M)$ is smooth and $\mathcal{N}(M)\supset \mathbf{C}^\ast=:T$. The torus $T$ acts with finitely many fixed points. It hence follows by Bialinicky--Birula results \cite{BB} (see also \cite[section~2.4]{Chriss}) that it has a cellular decomposition into attracting sets of its $T$--fixed points.
\end{proof}

In the rest of the paper we mainly concentrate on the quiver Grassmannians associated with indecomposable regular $Q$--representations.

\subsection{The variety $X=Gr_\mathbf{e}(R_n)$}\label{Sec:ActionOnX}

From now on we will focus on quiver Grassmannians associated with indecomposable regular $Q$--representations. It is not hard to show that $Gr_{\mathbf{e}}(R_n(\lambda))=Gr_{\mathbf{e}}(R_n(\mu))$ for every $\lambda,\mu\in k\cup\{\infty\}$ (see e.g. \cite{String}) and we hence consider the variety $X:=Gr_\mathbf{e}(R_n(0))$. It follows from the definition that
$$
X=\{N_1\subset N_2\subset k^n|\;J_n(0)N_1\subset N_2,\;\textrm{dim}N_i=e_i,\;i=1,2\}
$$
and hence $X$ is a closed subvariety of a partial flag variety. The group $\mathcal{N}:=\mathcal{N}(R_n)$ is the following:
$$
\mathcal{N}:=\{A\in GL_n(\mathbf{C})|\; AJA^{-1}=\lambda J,\;\textrm{for some }\lambda\in \mathbf{C}^\ast\}
$$
where $J:=J_n(0)$ and it acts on $X$ as $A\cdot(N_1,N_2)=(AN_1,AN_2)$.
\begin{lemma}\label{Lemma:N=UTZ}
The group $\mathcal{N}$ is the semi--direct product
$$
\mathcal{N}=\mathcal{U}\rtimes T\times Z
$$
where $\mathcal{U}$ is the unipotent radical of $\mathcal{N}$ of unipotent triangular Toeplitz matrices given by:
$$
\mathcal{U}:=\{\mathbf{1}_n+\sum_{i=1}^{n-1}a_iJ^i_n(0):\;a_i\in\mathbf{C}\},
$$
$T$ is the one--dimensional torus $T=\{t_\lambda:\;\lambda\in\mathbf{C}^\ast\}$ where $t_\lambda$ is the diagonal matrix
$$
t_\lambda:=\textrm{diag}(1,\lambda,\lambda^2,\cdots, \lambda^{n-1})
$$ and $Z$ consists of central elements $\{a_0\mathbf{1}_n|\;a_0\in\mathbf{C}^\ast\}$.
\end{lemma}
\begin{proof}
Every element of $\mathcal{U}TZ$ belongs to $\mathcal{N}$. Viceversa let $A\in\mathcal{N}$. Then it is easy to see that the columns $\mathbf{a}_1,\cdots,\mathbf{a}_n$ of $A$ satisfy the relation ${a}_{k+1}=J\mathbf{a}_{k}$ and hence $A\in \mathcal{U}TZ$. It is now easy to see that $\mathcal{U}$ is normal in $\mathcal{U}TZ$.
\end{proof}
For example for $n=5$, an element $A$ of the group $\mathcal{U}$ and a element $t_\lambda$ of $T$ have the form:
$$
\begin{array}{cc}
A=\left(
  \begin{array}{ccccc}
    1 & 0 & 0 & 0 & 0 \\
    a_1 & 1 & 0 & 0 & 0 \\
    a_2 & a_1 & 1 & 0 & 0 \\
    a_3 & a_2 & a_1 & 1 & 0 \\
    a_4 & a_3 & a_2 & a_1 & 1
  \end{array}
\right)&
t_\lambda=\left(
  \begin{array}{ccccc}
    1 & 0 & 0 &0 &0 \\
    0 & \lambda & 0 &0 &0 \\
    0 & 0 & \lambda^2 & 0&0 \\
    0 & 0&0 & \lambda^3 &0 \\
    0& 0 & 0&0 & \lambda^4 \\
  \end{array}
\right)
\end{array}
$$
for $a_1,a_2,a_3,a_4\in\mathbf{C}$ and $\lambda\in\mathbf{C}$. 
% As a consequence of lemma~\ref{Lemma:N=UTZ} the one--dimensional torus $T$ acts (algebraically) on $X$. We denote by $X^T$ the set of all its fixed points. In \cite{String}, \cite{CZ}, \cite{ZelElettronics} the cardinality of $X^T$, which is the Euler--Poincar\'e characterisitic $\chi(X)$ of $X$, is computed and it is given by:
% $$
% \chi(X)={e_2\choose e_1}{n-e_1\choose e_2-e_1}.
% $$
Clearly $L\in X^T$ if and only if $L$ is a coordinate sub--representation of $R_n$, i.e. both $L_1$ and $L_2$ are coordinate subspaces of $k^n$. In the next section we will encode this information in a combinatorial tool which is called the coefficient quiver of $R_n$.

We conclude this section by pointing out a useful isomorphism: 
\begin{equation}\label{Eq:GrRGrRAst}
\xymatrix@R=6pt{
 Gr_{\mathbf{e}}(R_n)\ar^{\varphi_n}[r]&Gr_{\mathbf{e}^\ast}(R_n)\\
N\ar@{|->}[r]&N^0
}
\end{equation}
where $(e_1,e_2)^\ast:=(n-e_2,n-e_1)$ and 
$
N^0:=\{f\in R_n^\ast|\;f(v)=0\;\forall v\in N\}
$
by using the identification $R_n\simeq R_n^\ast$. The torus $T$ acts on $R_n^\ast$ with contragredient action and this gives an action on $Gr_{\mathbf{e}^\ast}(R_n^\ast)$. Under the identification $R_n\simeq R_n^\ast$ these two actions differ by a character and so the identification $Gr_\mathbf{e}(R_n^\ast)\simeq Gr_\mathbf{e}(R_n)$ is $T$--equivariant. It is an easy check that isomorphism~\eqref{Eq:GrRGrRAst} is $T$--equivariant and involutive.

\subsection{The coefficient quiver of $R_n$}\label{Sec:CoefQuiver}
Following \cite{RingelExceptional} we associate to a $Q$--representation $M$ with standard basis $\mathcal{B}$ a quiver $\Gamma(M,\mathcal{B})$ called the coefficient quiver of $M$ in the basis $\mathcal{B}$. By definition $\Gamma(M,\mathcal{B})$ has the elements of $\mathcal{B}$ as vertices  and there is an arrow labeled by $a$ (resp. $b$) between two vertices $v$ and $v'$ if the coefficient of $v'$ in $m_a(v)$ (resp. $m_b(v)$) is non--zero. We call $\Gamma(M,\mathcal{B})$ the coefficient quiver of $M$ in the basis $\mathcal{B}$. The following are the coefficient quivers of some indecomposable $Q$--representations for $n=4$ in the standard basis. We denote by $k^{(i)}$ the vertex corresponding to $v_k^{(i)}$, $i=1,2$.
%\begin{figure}
$$
\Gamma(P_4)=
\xymatrix@C=6pt@R=8pt{
       &1^{(1)}\ar^a[dl]\ar_b[dr]&       &2^{(1)}\ar^a[dl]\ar_b[dr]&       &3^{(1)}\ar^a[dl]\ar_b[dr]&       &4^{(1)}\ar^a[dl]\ar_b[dr]&\\
1^{(2)}&       &2^{(2)}&       &3^{(2)}&       &4^{(2)}&       &5^{(2)}
}
$$
$$
\Gamma(R_5(0))=
\xymatrix@C=6pt@R=8pt{
       &1^{(1)}\ar^a[dl]\ar_b[dr]&       &2^{(1)}\ar^a[dl]\ar_b[dr]&       &3^{(1)}\ar^a[dl]\ar_b[dr]&       &4^{(1)}\ar^a[dl]\ar_b[dr]&&5^{(1)}\ar^a[dl]\\
1^{(2)}&       &2^{(2)}&       &3^{(2)}&       &4^{(2)}&       &5^{(2)}
}
$$
% $$
% \Gamma(R_5(\lambda))=
% \xymatrix@C=6pt@R=8pt{
%        &1^{(1)}\ar^a[dl]\ar@<1ex>^a[dr]\ar_b[dr]&       &2^{(1)}\ar^a[dl]\ar@<1ex>^a[dr]\ar_b[dr]&       &3^{(1)}\ar^a[dl]\ar@<1ex>^a[dr]\ar_b[dr]&       &4^{(1)}\ar^a[dl]\ar@<1ex>^a[dr]\ar_b[dr]&&5^{(1)}\ar^a[dl]\\
% 1^{(2)}&       &2^{(2)}&       &3^{(2)}&       &4^{(2)}&       &5^{(2)}
% }
% $$
$$
\Gamma(I_4)=
\xymatrix@C=6pt@R=8pt{
       &1^{(1)}\ar_b[dr]&       &2^{(1)}\ar^a[dl]\ar_b[dr]&       &3^{(1)}\ar^a[dl]\ar_b[dr]&       &4^{(1)}\ar^a[dl]\ar_b[dr]&&5^{(1)}\ar^a[dl]\\
&       &1^{(2)}&       &2^{(2)}&       &3^{(2)}&       &4^{(2)}
}
$$
%\end{figure}
In all these cases a one dimensional torus $T$ acts on quiver Grassmannians and the fixed points of this action are in bijection with successor closed subquivers of the corresponding coefficient quiver.

Let us consider $\Gamma(R_n)$. In lemma~\ref{Lemma:N=UTZ} we have seen the torus $T$ acts on $X$ by $\lambda\cdot v_k^{(i)}=\lambda^{k-1}v_k^{(i)}$ ($i=1,2$, $k\in[1,n]$). For every $r\geq1$ there exists a unique regular sub--representation of $R_n$ isomorphic to $R_r$ and it has the following property:
\begin{equation}\label{Def:Repulsive}
\textrm{if }N\in X\textrm{ is such that } \lim_{\lambda\rightarrow 0}t_\lambda N=R_r\textrm{ then } N=R_r.
\end{equation}
Indeed  this sub--representation is coordinate and lies in the extreme right hand side of $\Gamma(R_n)$. In particular the standard basis elements which generate $R_r$ have maximal weights. 

Roughly speaking the flow for $\lambda\rightarrow 0$ goes from right to left in $\Gamma(R_n)$. For example the line $\langle v_1^{(2)}+v_2^{(2)}\rangle$ generated by the vector $v_1^{(2)}+v_2^{(2)}$ goes to the line generated by $v_1^{(2)}$ as follows:
$$
\lim_{\lambda\rightarrow 0}\lambda\cdot(\langle v_1^{(2)}+v_2^{(2)}\rangle)=\lim_{\lambda\rightarrow 0}(\langle v_1^{(2)}+\lambda v_2^{(2)}\rangle)=\langle v_1^{(2)}\rangle.
$$

\subsection{Action of the torus on Hom--spaces}\label{Sec:HomSpaces}
Let $L$ and $L'$ be indecomposable $Q$--representations. The torus $T=\{t_\lambda:\;\lambda\in \mathbf{C}^\ast\}$ acts on the quiver Grassmannians associated with them. The action of $T$ naturally extends to the vector space $Hom(L,L')$ as follows: for $f\in Hom(L,L')$, $(t_\lambda f)(l):=t_\lambda f (t_{\lambda^{-1}}l)$. Following \cite{CrawleyTree} we endow the vector space $Hom(L,L')$ with a distinguished basis. Since $Hom(-,-)$ is additive, we assume that both $L$ and $L'$ are indecomposable. Let $\Gamma_L$ and $\Gamma_{L'}$ denote the coefficient quiver of $L$ and $L'$ respectively. We consider the set $(L,L')$ of triples $(\gamma,\alpha,\gamma')$ such that $\gamma$ is a connected \emph{predecessor closed} subquiver of $\Gamma_L$ (i.e. for every vertex $v\in\gamma$, every arrow $c: v'\rightarrow v$ with target $v$ belongs to $\gamma$) and $\gamma'$ is a successor closed subquiver of $\Gamma_{L'}$ and $\alpha :\gamma_0\rightarrow\gamma'_0$ is a bijective map from the set $\gamma_0$ of vertices of $\gamma$ to set $\gamma_0'$ of vertices of $\gamma'$ such that for every arrow $\xymatrix{v_k^{(1)}\ar^a[r]&v_k^{(2)}}$ (resp. $\xymatrix{v_k^{(1)}\ar^b[r]&v_{k+1}^{(2)}}$) of $\gamma$ there is an arrow $\xymatrix{\alpha(v_k^{(1)})\ar^a[r]&\alpha(v_k^{(2)})}$ (resp. $\xymatrix{\alpha(v_k^{(1)})\ar^b[r]&\alpha(v_{k+1}^{(2)})}$) of $\gamma'$ with the same label. The next proposition is a special case of \cite{CrawleyTree} and we hence avoid to repeat the proof (which is quite simple in this case).

\begin{figure}
$$
\xymatrix@C=10pt@R=15pt{&&&&&&&&&&&&\\
&1\circ\ar_a[dl]\ar^b[dr]\ar@/^1pc/@{..>}[rrrrrrrrr]&&2\circ\ar_a[dl]\ar^b[dr]\ar@/^1pc/@{..>}[rrrrrrrrr]&&3\ar_a[dl]&&&1\ar_a[dl]\ar^b[dr]&&2\bullet\ar_a[dl]\ar^b[dr]&&3\bullet\ar_a[dl]\\
1\circ\ar@/_1pc/@{..>}[rrrrrrrrr]&&2\circ\ar@/_1pc/@{..>}[rrrrrrrrr]&&3&&&1&&2\bullet&&3\bullet&
}
$$
\caption{An element of the standard basis of $Hom(R_3,R_3)$: the circles (resp. bullets) highlight a predecessor (successor) closed subquiver $\gamma$ (resp. $\gamma'$) 
of $\Gamma(R_3)$. The dotted arrows 
show the corresponding $f_{\gamma\gamma'}$
}
\label{Fig:Hom}
\end{figure}
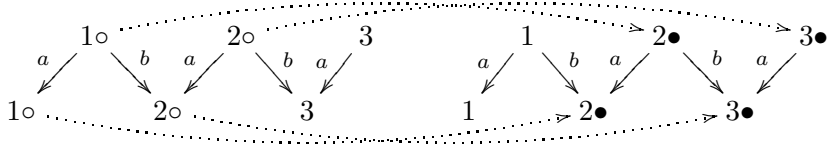

\begin{proposition}\label{Prop:BasisHom}
We consider the map\
$$
B:(L,L')\rightarrow Hom(L,L'):\;\;\; (\gamma,\alpha,\gamma')\mapsto f_{\gamma\gamma'}
$$
which associates to $(\gamma,\alpha,\gamma')$ the homomorphism
$$
f_{\gamma\gamma'}(v)=\left\{
\begin{array}{cc}
\alpha(v)&\textrm{ if }v\in \gamma_0,\\
0&\textrm{ otherwise. }
\end{array}\right.
$$
The image of $B$ is a basis of $Hom(L,L')$ which we call standard.
\end{proposition}
Figure~\ref{Fig:Hom} illustrates proposition~\ref{Prop:BasisHom}. The torus $T$ acts diagonally on the elements of the standard basis of $Hom(L,L')$ as follows: if $(\gamma,\alpha,\gamma')\in(L,L')$ and the vertices of $\gamma$ have consecutive weights $k$, $k+1$, $\cdots$ and the vertices of $\gamma'$ have consecutive weights $k'$, $k'+1$, $\cdots$ then
$$
t_\lambda f_{\gamma\gamma'}=\lambda ^{k-k'}f_{\gamma\gamma'}
$$
and we say that $f_{\gamma\gamma'}$ has weight $k-k'$. We say that an element $f_{\gamma\gamma'}$ has \emph{positive weight} if $k-k'>0$ and it has \emph{negative weight} if $k-k'<0$. We denote by
$Hom(L,L')^+$ the vector subspace of $Hom(L,L')$ spanned by standard basis elements with positive weight.
As an application of proposition~\ref{Prop:BasisHom} we compute the dimension of the Hom--spaces between indecomposable $Q$--representations. It is known (and not difficult to prove by using proposition \ref{Prop:BasisHom}) that $Hom(R_s,P_l)=Hom(I_s,P_l)=Hom(I_s,R_l)=0$ for all $s,l\geq0$. We hence consider the remaining cases.
\begin{lemma}\label{Lemma:Hom}
For every $l,s\geq0$ we have:
\begin{eqnarray}\label{Eq:Dim(P,P)}
\textrm{dim} Hom(P_s,P_l)&=& [l-s+1]_+\\
\label{Eq:Dim(P,R)}
\textrm{dim} Hom(P_s,R_l)&=& l\\
\label{Eq:Dim(P,I)}
\textrm{dim} Hom(P_s,I_l)&=& l+s\\
\label{Eq:Dim(R,R)}
\textrm{dim} Hom(R_s,R_l)&=& \textrm{min}(s,l)\\
\label{Eq:Dim(R,I)}
\textrm{dim} Hom(R_s,I_l)&=& s\\
\label{Eq:Dim(I,I)}
\textrm{dim} Hom(I_s,I_l)&=& [s-l+1]_+
\end{eqnarray}
where $[b]_+:=\textrm{max}(b,0)$.
\end{lemma}

We conclude this section by pointing out that the action of the torus $T$ on $Hom(L,L')$ affords an action of the torus on the space $Ext^1(L,L')$. With respect to this action long exact sequences in cohomology are $T$--equivariant. 
% 
% \begin{lemma}
% Let $L$ and $M$ be $Q$--representations such that a one--dimensional torus $T$ acts on $Hom(L,M)$. If $Hom(L,M)^+=0$ then $Ext^1(L,M)^+=0$.
% \end{lemma}
% \begin{proof}
% We choose a $T$--equivariant injective resolution of $M$
% \end{proof}

\subsection{Stratification of $X$}\label{Sec:Stratification}

Every sub--representation $N$ of $R_n$ is of the form $N=P\oplus R_r$ where $P$ is  preprojective and $R_r$ is an indecomposable (possibly zero) regular $Q$--representation, for some $r\geq0$. Similarly the quotient $R_n/N=R_{r'}\oplus I$ where $R_{r'}$ is regular indecomposable (possibly zero) and $I$ is preinjective. We hence can give the following definition.
\begin{definition}\label{Def:KN}
Let $N\in X$ with $N=P\oplus R_r$ and $R_n/N= R_{r'}\oplus I$ with $P$ preprojective, $I$ preinjective and some $r,r'\geq0$. We define the integer:
$$
K_N=K_N(X):=\textrm{min}(r,r')
$$
\end{definition}
It is easy to see that $K_N=\textrm{dim}Ext^1(N,R_n/N)$. Indeed
$$
\textrm{dim}Ext^1(N,R_n/N)=\textrm{dim}Ext^1(R_r,R_{r'})=\textrm{dim}Hom(R_r,R_{r'})=\textrm{min}(r,r')
$$ 
where in the first equality we use the well--known fact that $Ext^1(R,I)=Ext^1(P,R)=Ext^1(P,I)=0$ for every preprojective $P$, regular $R$ and preinjective $I$ $Q$--representations; in the second equality we use the AR--formula (see e.g. \cite{ASS});
% $$
% Ext^1(N,M)\simeq DHom(M,\tau N);
% $$
in the last equality we use \eqref{Eq:Dim(R,R)}. In particular it is known that
$$
\langle\mathbf{dim}N,\mathbf{dim}N'\rangle=\textrm{dim} Hom(N,N')-\textrm{dim}Ext^1(N,N')
$$
and in view of \eqref{Eq:TN} we get
\begin{equation}\label{Eq:DimTN}
\textrm{dim}T_N=\langle\mathbf{e},n\delta-\mathbf{e}\rangle+K_N.
\end{equation}
Equation \eqref{Eq:DimTN} implies that a point $N\in X$ is smooth if and only if either $N$ or $M_n/N$ do not have a regular direct summand.

The next theorem provides a stratification of $X$ and it is essential for our proof of the existence of a cellular decomposition of $X$. Let us define the strata.

\begin{definition}\label{Def:Xk}
For every integer $k\geq0$ define the set:
$$
X_k=X_k(X):=\{N\in X|\;K_N\geq k\}
$$
where $K_N$ is given in definition~\ref{Def:KN}.
\end{definition}

\begin{theorem}\label{Thm:Stratification}
\begin{enumerate}
\item The set $X_{k+1}$ is a closed $T$--stable subvariety of $X_k$. Moreover there is a $T$--equivariant isomorphism:
    \begin{equation}\label{Eq:XkGr}
    X_k\simeq Gr_{(e_1-k,e_2-k)}(R_{n-2k});
    \end{equation}
\item The sets $X_k$ provide a stratification of $X$
\begin{equation}\label{Eq:StratificationThm}
X=X_0\supset X_1\supset\cdots\supset X_s
\end{equation}
where $s=\textrm{min}(e_1,n-e_2)$;
\item The variety $X_k\setminus X_{k+1}$ is smooth (inside $X_k$) and
$$
\forall N\in X_k\setminus X_{k+1},\;\;\; \lim_{\lambda\rightarrow0}t_\lambda N\in X_k\setminus X_{k+1}
$$
\end{enumerate}
\end{theorem}
\begin{proof}
We consider the two subvarieties of $X$:
for $k\in[0,e_1]$
\begin{equation}\label{Eq:DefXk'}
X'_k(X):=\{N\in Gr_\mathbf{e}(R_n)|\; R_k\textrm{ is a subrepresentation of } N\}
\end{equation}
and for $k\in[0,n-e_2]$
\begin{equation}\label{Eq:DefXk''}
X''_k(X):=\{N\in Gr_\mathbf{e}(R_n)|\; R_k\textrm{ is a quotient of } R_n/N\}.
\end{equation}
It follows from the definitions that:
\begin{equation}\label{Eq:XkXk'Xk''}
X_k=X'_k\cap X''_k. 
\end{equation}
We now collect some properties of $X_k'$ and $X''_k$.
\begin{lemma}\label{Lemma:X'X''}
\begin{enumerate}
\item The isomorphism $\varphi_n $ defined in \eqref{Eq:GrRGrRAst} induces an isomorphism 
$$
X'_k(Gr_\mathbf{e}(R_n))\simeq X''_k(Gr_{\mathbf{e}^\ast}(R_n))
$$
and hence also
$$
X''_k(Gr_\mathbf{e}(R_n))\simeq X'_k(Gr_{\mathbf{e}^\ast}(R_n)).
$$
\item For every $k\in[0,e_1]$, $X'_k(X)$ is a $T$--stable subvariety of $X$. For every $N\in X'_k(X)$, 
$\lim_{\lambda\rightarrow0}t_\lambda N\in X'_k(X)$. There is a $T$--equivariant isomorphism 
\begin{equation}\label{Eq:Xk'Gr}
X'_k(X)\simeq Gr_{\mathbf{e}-k\delta}(R_{n-k})
\end{equation}

\item For every $k\in[0,n-e_2]$,  $X''_k(X)$ is a $T$--stable subvariety of $X$. For every $N\in X''_k(X)$, 
$\lim_{\lambda\rightarrow0}t_\lambda N\in X''_k(X)$. There is a $T$--equivariant isomorphism 
\begin{equation}\label{Eq:Xk''Gr}
X''_k(X)\simeq Gr_{\mathbf{e}}(R_{n-k})
\end{equation}
\end{enumerate}
\end{lemma}
\begin{proof}
Part $(1)$ follows by straightforward check. Part $(3)$ follows from part $(1)$ by using the isomorphism~\eqref{Eq:GrRGrRAst}. Part $(2)$ follows from property \eqref{Def:Repulsive} of $R_k$. It remains to check \eqref{Eq:Xk'Gr}. We consider the map:
$$
\xymatrix@R=6pt{
\{N\in Gr_\mathbf{e}(R_n)|\; R_k\leq N\}\ar[r]&Gr_{\mathbf{e}-k\delta}(R_{n-k})\\
N\ar@{|->}[r]&N/R_k
}
$$
Since $R_n$ is uniserial, there is a canonical short exact sequence:
\begin{equation}\label{Eq:SESRk}
\xymatrix{0\ar[r]& R_k\ar^{\iota_k}[r]& R_n\ar^{\pi_{n-k}}[r]&R_{n-k}\ar[r]&0}
\end{equation}
and hence the map is well defined and bijective. 
\end{proof}
We now conclude the proof of the theorem. Everything follows from lemma~\ref{Lemma:X'X''} by \eqref{Eq:XkXk'Xk''} except the smoothness of $X_k\setminus X_{k+1}$. To prove smoothness notice that it is sufficient to prove smoothness of $X_0\setminus X_1$ by \eqref{Eq:XkGr}, since $K_N(X_k)=K_N(X)-k$. But by \eqref{Eq:DimTN} $X_0\setminus X_1$ is the smooth locus of $X$.
\end{proof}
\begin{corollary}
The variety $X=Gr_\mathbf{e}(R_n)$ is smooth if and only if $e_1=0$ (in which case $X\simeq Gr_{e_2}(k^n)$) or $e_2=n$ (in which case $X\simeq Gr_{e_1}(k^n)$).
\end{corollary}
\begin{proof}
By theorem~\ref{Thm:Stratification}, $X$ is smooth if and only if $X=X_0$ if and only if $s=\textrm{min}(e_1,n-e_2)=0$.
\end{proof}
\begin{corollary}
 The quiver Grassmannian $Gr_{(e_1,e_2)}$ has dimension:
$$
\textrm{dim}Gr_{(e_1,e_2)}=\langle \mathbf{e}, n\delta-\mathbf{e}\rangle=(e_2-e_1)(n-(e_2-e_1))
$$
\end{corollary}
\begin{proof}
 It follows from \eqref{Eq:DimTN}.
\end{proof}

\subsection{Cellular decomposition of $X$}\label{Sec:CellDec}
In this section we provide a cellular decomoposition of $X=Gr_\mathbf{e}(R_n)$. Following \cite{BB} (see also \cite{Chriss}), for every fixed point $L\in X^T$ we consider its \emph{attracting set} defined as follows:
\begin{equation}\
X_L:=\{N\in X|\; \lim_{\lambda\rightarrow0} t_\lambda N=L\}.
\end{equation}
In particular $L\in X_L$ for every $L\in X^T$.
\begin{theorem}\label{Thm:CellDec}
For every $L\in X^T$ the corresponding attracting set $X_L$is an affine space and $X_L\simeq T_L(X_L)\simeq T_L^+(X):=Hom(L,R_n/L)^+$. Moreover
$$
X=\cup_{L\in X^T} X_L
$$
\end{theorem}
\begin{proof}
Bialynicky--Birula's results on cellular decomposition of a projective variety $X$ \cite{BB} continues to hold if the variety is smooth but only quasi--projective provided that the action of the torus is such that $\lim_{\lambda\rightarrow 0}t_\lambda N$ belongs to $X$ for every $N\in X$. This is an easy consequence of Hironaka's resolution of singularities. Since of theorem~\ref{Thm:Stratification} we apply this to $X_0\setminus X_1$ and we get the result.
\end{proof}
\subsection{Description of the cells}\label{Sec:Cells}
In this section we describe the cell $X_L$ associated with every $L\in X^T$ (see theorem~\ref{Thm:CellDec}). If $L\simeq P_{r}$ is an indecomposable preprojective subrepresentation of $R_n$ of dimension $(r,r+1)$, we write $L:={}_k(P_{e_1})$ if $L$ is generated by $v_k^{(1)}, v_{k+1}^{(1)},\cdots,v_{k+r-1}^{(1)}$ and $v_k^{(2)}, v_{k+1}^{(2)},\cdots, v_{k+r}^{(2)}$. For example the following figure shows the subrepresentation ${}_2(P_1)$ of $R_5$:
$$
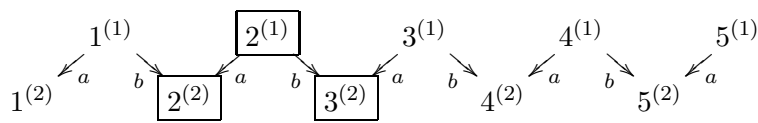
\xymatrix@C=6pt@R=8pt{
&1^{(1)}\ar^a[dl]\ar_b[dr]&       &*+[F]{2^{(1)}}\ar^a[dl]\ar_b[dr]&       &3^{(1)}\ar^a[dl]\ar_b[dr]&       &4^{(1)}\ar^a[dl]\ar_b[dr]&&5^{(1)}\ar^a[dl]\\
1^{(2)}&       &*+[F]{2^{(2)}}&       &*+[F]{3^{(2)}}&       &4^{(2)}&       &5^{(2)}
}
$$
\begin{theorem}\label{Thm:Cells}
\begin{enumerate}
\item If $L\in X^T$ is indecomposable then $X_L=\mathcal{U}L$. In particular if $L={}_k(P_{e_1})$ then $\textrm{dim}X_L=n-k$;
\item If $L=L'\oplus L''$ with $Hom(L',L'')^+=Hom(L',L'')$ then
\begin{equation}\label{Eq:DimXL}
\textrm{dim} X_L=\textrm{dim} X_{L'}+\textrm{dim} X_{L''}-\langle\mathbf{dim}L',\mathbf{dim}L''\rangle
\end{equation}
\end{enumerate}
\end{theorem}
\begin{proof}
We prove part (1). The unipotent group $\mathcal{U}$ defined in lemma~\ref{Lemma:N=UTZ} is a subgroup of dimension $n-1$ of the group of $n\times n$ unipotent triangular matrices. Let $L$ be an indecomposable
sub--representation of $R_n$ of dimension vector $\mathbf{e}$. If $L=R_r$ then $\mathbf{e}=(r,r)$ and $L_1=L_2$ is the vector subspace of $k^n$ spanned by the last $r$ standard basis vectors. In particular $AL=L$ for every $A\in \mathcal{U}$. On the other hand we have already noticed that regular sub--representations have property \eqref{Def:Repulsive} and hence $X_L=\{L\}$ and we get $X_L=\mathcal{U}L$ if $L$ is regular. Let $L$ be an indecomposable preprojective sub--representation of $R_n$. Then $e_2=e_1+1$ and we assume that $L_1$ is spanned by consecutive standard basis vectors $\{v_k, v_{k+1},\cdots, v_{k+e_1-1}\}$ and $L_2$ is spanned by $L_1$ and $v_{k+e_1}$. The following easy lemma says that $\textrm{dim}\, \mathcal{U}L=n-k$.
\begin{lemma}
Let $W$ be a vector subspace of $k^n$ spanned by some basis vectors $\{v_i:\; i\geq k\}$ and $v_k$ is in this set. Then the stabilizer $Stab_{\mathcal{U}}(W)=Stab_\mathcal{U}(<v_k>)$. In particular the orbit $\mathcal{U}W$ has dimension $n-k$.
\end{lemma}
\begin{proof}
We prove that $Stab_{\mathcal{U}}(W)\supset Stab_\mathcal{U}(<v_k>)$ being the other inclusion clear.
The columns $\mathbf{a}_1,\cdots, \mathbf{a}_n$ of a matrix $A\in\mathcal{U}$ satisfy the relation: $\mathbf{a}_l=J_n(0)\mathbf{a}_{l-1}$ for every $l=2,3,\cdots, n$. Let $A\in Stab_\mathcal{U}(<v_k>)$ and $A v_k=\beta v_k$ for some non--zero scalar $\beta$. Let $v:=\sum _{i\geq 0}b_i v_{k+i}$ be an element of $W$; we get:
$$
Av=\sum_{i\geq 0}b_i Av_{k+i}=\sum_{i\geq 0}b_i \mathbf{a}_{k+i}=\sum_{i\geq 0}b_i J^i\mathbf{a}_{k}=\sum_{i\geq0}b_i J^i A v_k=\sum_i b_i\beta J^i v_k
$$
which belongs to $W$ and hence $A\in Stab_{\mathcal{U}}(W)$. In particular $\mathcal{U}W=\mathcal{U}(<v_k>)$ has dimension $n-k$.
\end{proof}
We now prove that $\textrm{dim} X_L=n-k$ and since $\mathcal{U}L\subset X_L$ we get the equality. By theorem~\ref{Thm:CellDec} we have to compute $\textrm{dim} T_L^+(X)=\textrm{dim} Hom(L,R_n/L)^+$. Since $L$ is indecomposable  the quotient $R_n/L$ is the direct sum of at most two indecomposables as follows
$$
R_n/L=R_{k-1}\oplus I_t
$$
where $t=n-1-k-e_1$. We have:
$$
Hom(L,R_n/L)^+=Hom(P_{e_1}, I_t).
$$
In view of \eqref{Eq:Dim(P,I)}, $\textrm{dim}(Hom(L,R_n/L)^+)=e_1+t=n-k$ and we are done.

We now prove part (2). There are short exact sequences:
\begin{equation*}
\xymatrix{0\ar[r]& L'\ar[r]& R_n/L''\ar[r]&R_{n}/(L'\oplus L'')\ar[r]&0}
\end{equation*}
and
\begin{equation*}
\xymatrix{0\ar[r]& L''\ar[r]& R_n/L'\ar[r]&R_{n}/(L'\oplus L'')\ar[r]&0}
\end{equation*}
The torus acts on Hom--spaces and on Ext--spaces between $T$--fixed points as explained in section~\ref{Sec:HomSpaces}. Moreover these spaces split into the subspaces spanned by standard basis vectors with positive and negative weights respectively. We apply the functors $Hom(L'',-)$ and $Hom(L',-)$ to the previous short exact sequences, we take the positive part and we get the following exact sequences:
\begin{equation*}
\xymatrix{0\ar[r]& 0\ar[r]& T_{L''}(X_{L''})\ar[r]&Hom(L'',R_n/L)\ar[r]&0}
\end{equation*}
and
\begin{equation*}
\xymatrix@C=9pt{0\ar[r]& Hom(L',L'')\ar[r]& T_{L'}(X_{L'})\ar[r]&Hom(L',R_n/L)\ar[r]&Ext^1(L',L'')^+\ar[r]&0}
\end{equation*}
where we use the following equalities:
\begin{equation}
Ext^1(L'',L')^+=0=Ext^1(L',R_n/L')^+
\end{equation}
\begin{equation}
Ext^1(L',L'')^+=Ext^1(L',L'')
\end{equation}
% \begin{equation}
% Hom(L',L'')^+=Hom(L',L'')
% \end{equation}
This follows by the fact that if  $Hom(M,L)^+=0$ then $Ext^1(M,L)^+=0$; indeed one can always take a minimal $T$--equivariant injective resolution of $L$ and apply the functor $Hom(M,-)$. 
By summing up we get the short exact sequence:
\begin{equation*}
\xymatrix@C=9pt{0\ar[r]& Hom(L',L'')\ar[r]& T_{L'}(X_{L'})\oplus T_{L''}(X_{L''})\ar[r]&T_L(X_L)\ar[r]&Ext^1(L',L'')\ar[r]&0}
\end{equation*}
Since for every fixed point $L$, $\textrm{dim} X_L=\textrm{dim} T_L(X_L)$ we get \eqref{Eq:DimXL}.
\end{proof}

\subsection{Betti numbers}\label{Sec:Betti}

We now use the results of the previous sections in order to compute the Betti numbers of $X=Gr_\mathbf{e}(R_n)$. Since $X$ has a cellular decomposition (theorem~\ref{Thm:CellDec}) the odd cohomology spaces are zero and the $2i$--th Betti number $b_{2i}=b_{2i}(X):=\textrm{dim}H^{2i}(X)$ equals the number of cells of dimension $i$. 

Before stating the main result of this section we start with the special case $e_2=e_1+1$. 
\begin{theorem}
Let $X=Gr_{(e_1,e_1+1)}(R_n)$. The even Betti numbers of $X$ are the following:
\begin{equation}\label{Eq:Betti}
b_{2i}=\left\{
\begin{array}{cc}
i+1&\textrm{ if }0\leq i\leq s,\\
s+1&\textrm{ if }s\leq i\leq n-1-s,\\
n-i&\textrm{ if }n-1-s\leq i\leq n-1\\
\end{array}
\right.
\end{equation}
where $s=min(e_1,n-1-e_1)$.

The Poincar\'e polynomial $P_X(t):=\sum_{i=0}^{\textrm{dim}X}b_it^i$ of $X$ equals:
\begin{equation}\label{Eq:PxSpecial}
P_X(t^{1/2})=\left(\frac{t^{e_1+1}-1}{t-1}\right)\left(\frac{t^{n-e_1}-1}{t-1}\right)
\end{equation}
\end{theorem}
\begin{proof}
For $k\in[0,e_1]$ let us consider the variety $X'_k$ defined in \eqref{Eq:DefXk'} and the difference 
\begin{equation}\label{Eq:Y'k}
Y'_k=Y'_k(X):=X'_k(X)\setminus X'_{k+1}(X).
\end{equation}
By definition $N\in Y_k$ if $N=R_k\oplus P$ for some preprojective $P$. By lemma~\ref{Lemma:X'X''} we have
\begin{equation}\label{Eq:HXHXk}
\textrm{dim}H^i(X)=\sum _{k=0}^{e_1} \textrm{dim}H^i_c(Y'_k)
\end{equation}
where $H_c^i(Y'_k)$ is the $i$--th cohomology space of $Y'_k$ with compact support. 
Moreover $\textrm{dim}H^i_c(Y'_k)=0$ for $i$ odd and
\begin{equation}\label{Eq:BettiX'}
\textrm{dim}H^{2i}_c(Y'_k)=\left\{
\begin{array}{cc}
1&\textrm{ if }i\in[e_1-k,n-k-1]\\
0&\textrm{ otherwise.}
\end{array}
\right. 
\end{equation}
and hence \eqref{Eq:Betti} follows from \eqref{Eq:HXHXk}. It remains to prove \eqref{Eq:BettiX'}.
It is easy to see that 
\begin{equation}
Y'_k(X)\simeq Y'_0(Gr_{e_1-k,e_2-k}(R_{n-k})).
\end{equation}
% The number $\textrm{dim}H^{2i}_c(Y'_k)$ equals the number of attracting sets $X_L$, $L\in X'_k\cap X^T$, of dimension $i$. 
The elements of $Y'_0(Gr_{(e_1,e_1+1)}(R_n))$ are all the indecomposable preprojective subrepresentations of $R_n$ of 
dimension 
$(e_1,e_1+1)$ and hence they are all isomorphic to $P_{e_1}$. We denote by ${}_i(P_{e_1})$ the fixed point of 
$Y'_0(Gr_{(e_1,e_1+1)}(R_n))$ generated by the vectors $\{v_i^{(1)}, v_{i+1}^{(1)},\cdots,v_{i+e_1-1}^{(1)}\}$ and 
$\{v_i^{(2)}, v_{i+1}^{(2)},\cdots,v_{i+e_1}^{(2)}\}$ for $i\in[1,n-e_1]$. By theorem~\ref{Thm:Cells} 
$\textrm{dim}X_{{}_i(P_{e_1})}=n-i\in[e_1,n-1]$ 
and hence for every $i\in[e_1,n-1]$ there is a unique cell of dimension $i$ and we get \eqref{Eq:BettiX'}. 
\end{proof}
Notice that \eqref{Eq:PxSpecial} can be written as follows:
$$
P_X(t)=P_{Gr_{1}(e_2)}(t)P_{Gr_{1}(n-e_1)}(t)
$$
where $Gr_s(t)$ denotes the Grassmannians of $s$--dimensional vector subspaces of a $t$--dimensional vector space. Surprisingly this turns out to be a general fact and it is the main result of this section. 

\begin{theorem}\label{Thm:Poincare}
Let $X=Gr_{(e_1,e_2)}(R_n)$. The Poincar\'e polynomial $P_X(t)$ of $X$ equals
\begin{equation}\label{Eq:MainThmBetti}
P_X(t)=P_{Gr_{(e_2-e_1)}(e_2)}(t)P_{Gr_{(e_2-e_1)}(n-e_1)}(t).
\end{equation}
\end{theorem}
\begin{proof}
We proceed by induction on $n\geq e_1\geq0$. For $e_1=0$, $X\simeq Gr_{e_2}(n)$ and \eqref{Eq:MainThmBetti} follows. Let $1\leq e_1\leq e_2$. The variety $X$ can be decomposed in 
$$
X=Y'_0(X)\cup X'_1(X)
$$
where $Y'_0=Y'_0(X)$ consists of all preprojective subrepresentations of $R_n$ in $X$ and $X'_1$ consists of subrepresentations of $R_n$ in $X$ having a non--zero regular subrepresentation (see \eqref{Eq:DefXk'} and \eqref{Eq:Y'k}). Moreover, by \eqref{Eq:Xk'Gr}, $X'_1(X)\simeq Gr_{(e_1-1,e_2-1)}(R_{n-1})$ and we hence have:
\begin{equation}\label{Eq:Px=Py}
 P_X(t)=P_{Y'_0}(t)+P_{Gr_{(e_1-1,e_2-1)}(R_{n-1})}(t).
\end{equation}
By inductive hypothesis the equality \eqref{Eq:MainThmBetti} holds if and only if the following equality holds: 
\begin{equation}\label{Eq:Py=Px}
P_{Y'_0}(t)=\left(P_{Gr_{(e_2-e_1)}(e_2)}(t)-P_{Gr_{(e_2-e_1)}(e_2-1)}(t)\right)P_{Gr_{(e_2-e_1)}(n-e_1)}(t).
\end{equation}
We hence prove \eqref{Eq:Py=Px}. We make the following choice: we fix a basis $v_1\cdots v_{s}$ of an $s$--dimensional vector space and we let the torus act on $Gr_t(s)$ by $t_\lambda v_i=\lambda^i v_i$. We consider the vector subspace of $k^{e_2}$ generated by $v_2,\cdots, v_{e_2}$ and the corresponding embedding $Gr_{(e_2-e_1)}(e_2-1)\subset Gr_{(e_2-e_1)}(e_2)$. With this choices the difference $Gr_{(e_2-e_1)}(e_2)\setminus Gr_{(e_2-e_1)}(e_2-1)$ is $T$--stable and for every  point $W$ of it, $\lim_{\lambda\to 0}t_\lambda W$ still belongs to it. The right hand side of \eqref{Eq:Py=Px} is the Poincar\'e polynomial (with respect to the cohomology with compact support) of the smooth projective variety:
$$
G:=\left(Gr_{(e_2-e_1)}(e_2)\setminus Gr_{(e_2-e_1)}(e_2-1)\right)\times Gr_{(e_2-e_1)}(n-e_1)
$$
The one--dimensional torus $T$ acts on $G$ and the attracting sets  of the $T$--fixed points form a cellular decomposition of $G$. We prove that there exists a bijection between the cells of $Y'_0$ of dimension $k$ and the cells of $G$ of dimension $k$. 

A point of $Y'_0$ is a direct sum of precisely $(e_2-e_1)$ preprojective subrepresentations of $R_n$ (this follows by considering their dimension vector). The $T$--fixed points have the form:
$$
L:={}_{k_1}(P_{r_1})\oplus {}_{k_2}(P_{r_2})\oplus\cdots\oplus {}_{k_{(e_2-e_1)}}(P_{r_{(e_2-e_1)}})
$$ 
where $r_1+\cdots +r_{e_2-e_1}=e_1$, $r_i\geq0$ and ${}_k(P_r)$ denotes the unique indecomposable preprojective subrepresentation of $R_n$ of dimension vector $(r,r+1)$ generated by $v_k^{(1)}, v_{k+1}^{(1)},\cdots,v_{k+r-1}^{(1)}$ and $v_k^{(2)}, v_{k+1}^{(2)},\cdots, v_{k+r}^{(2)}$. 
In view of theorem~\ref{Thm:Cells} it is easy to see that the dimension of the attracting cell of $L$ equals
$$
\textrm{dim}X_L=n(e_2-e_1)-\sum_{i=1}^{e_2-e_1}k_i-(e_2-e_1)^2+\sum_{i=1}^{e_2-e_1}\sum_{j=1}^{i}(r_j-r_i+1).
$$
We consider the set $\alpha(\mathbf{e},k,n)$ which parametrezies the $T$--fixed points of $Y'_0$ whose attracting set has dimension $k$, i.e.
$$
\begin{array}{c}
\alpha(\mathbf{e},k,n):=\{(k_1,k_2,\cdots,k_{e_2-e_1},r_1,r_2,\cdots,r_{e_2-e_1})|\;\\
1\leq k_1\leq k_1+r_1< k_2\leq k_2+r_2\langle\cdots<k_{e_2-e_1}\leq k_{e_2-e_1}+r_{e_2-e_1}\leq n;\\
r_1+r_2+\cdots r_{e_2-e_1}=e_1,\; r_i\geq0;\\
n(e_2-e_1)-\sum_{i=1}^{e_2-e_1}k_i-(e_2-e_1)^2+\sum_{i=1}^{e_2-e_1}\sum_{j=1}^{i}(r_j-r_i+1)=k
\}.
\end{array}
$$

On the other hand we consider the $T$--fixed points of $G$ and their attracting sets. The $T$--fixed points of $Gr_t(s)$ 
consist of coordinate vector subspaces of dimension $t$ and they are naturally parametrized by tuples $(a_1,\cdots a_t)$ of integers $1\leq a_1<\cdots< a_t\leq s$. The corresponding cell $\mathcal{O}_{(a_1,\cdots, a_t)}$ has dimension 
$$
s-a_1-(t-1)+s-a_2-(t-2)+\cdots s-a_t.=ts-\sum_{i=1}^{t}a_i-\sum_{i=1}^{t-1}i.
$$ 
The $T$--fixed points of $Gr_{(e_2-e_1)}(e_2)\setminus Gr_{(e_2-e_1)}(e_2-1)$ are the coordinate vector subspaces of $k^{e_2}$  containing $v_1$. The following set hence parametrizes the cells of $G$ of dimension $k$:
$$
\begin{array}{c}
\beta(\mathbf{e},k,n):=\{(a_1,a_2,\cdots,a_{e_2-e_1},b_2,b_3,\cdots,b_{e_2-e_1})|\;\\
1\leq a_1< a_2\langle\cdots<a_{e_2-e_1}\leq n-e_1;\\
2\leq b_2< \cdots<b_{e_2-e_1}\leq e_2;\\
n(e_2-e_1)-\sum_{i=1}^{e_2-e_1}a_i-\sum_{i=2}^{e_2-e_1}b_i+(e_2-e_1-1)=k
\}.
\end{array}
$$
We consider the map 
\begin{eqnarray}
a_1&=&k_1\\
a_i&=& k_i-r_1-r_2-\cdots -r_{i-1}\\
b_i&=&\sum_{j=0}^{i-2}r_{(e_2-e_1-j)}+i
\end{eqnarray}
for $i\in[2,e_2-e_1]$. It is straightforward to verify that this map is a bijection between $\alpha(\mathbf{e},k,n)$ and $\beta(\mathbf{e},k,n)$. It follows that $Y'_0$ and $G$ have the same Betti numbers and hence \eqref{Eq:Py=Px} follows.
\end{proof}

\begin{corollary}
The Poincar\'e polynomial of a quiver Grassmannian associated with the indecomposable preprojective $P_n$ and the indecomposable preinjective $I_n$ ($n\geq0$) are respectively the following:
\begin{equation}\label{Eq:PoiPn}
P_{Gr_\mathbf{e}(P_n)}(t)=P_{Gr_{e_1}(e_2-1)}(t)P_{Gr_{(e_2-e_1)}(n+1-e_1)}(t)
\end{equation}
\begin{equation}\label{Eq:PoiIn}
P_{Gr_\mathbf{e}(I_n)}(t)=P_{Gr_{e_1}(e_2+1)}(t)P_{Gr_{(e_2-e_1)}(n-e_1)}(t)
\end{equation}
\end{corollary}
\begin{proof}
The equality \eqref{Eq:PoiIn} follows from \eqref{Eq:PoiPn} by the isomorphism:
$$
Gr_{(e_1,e_2)}(I_n)\simeq Gr_{(n-e_2,n+1-e_1)}(P_n).
$$
We hence prove \eqref{Eq:PoiPn}. As in the proof of theorem~\ref{Thm:Poincare}, let $Y'_0=Y'_0(Gr_{\mathbf{e}}(R_{n+1}))$ be the subvariety of $Gr_{\mathbf{e}}(R_{n+1})$ of all preprojective subrepresentations and let 
$$
G':=Gr_{(e_2-e_1)}(e_2)\setminus Gr_{(e_2-e_1)}(e_2-1)
$$
with the convention that $Gr_{(e_2-e_1)}(e_2-1)$ consists of all the elements of $Gr_{(e_2-e_1)}(e_2)$ not containing the first basis vector (as in the proof of theorem~\ref{Thm:Poincare}). In view of \eqref{Eq:Py=Px}  it is sufficient to prove the following equalities:
\begin{eqnarray}\label{Eq:PG'e}
 P_{G'}(t)&=&t^{2e_1}P_{Gr_{(e_2-e_1-1)}(e_2-1)}(t),\\\label{Eq:PY'e}
P_{Y'_0}(t)&=&t^{2e_1}P_{Gr_{\mathbf{e}}(P_{n})}(t).
\end{eqnarray}
The proof of \eqref{Eq:PG'e} is similar to the proof of theorem~\ref{Thm:Poincare}: there is an obvious bijection between the cells of $G'$ of dimension $k$ and the cells of $Gr_{(e_2-e_1)}(e_2-1)$ of dimension $k-e_1$. Let us prove \eqref{Eq:PY'e}. Let $L\in Y'_0$. Then $L$ is a sum of preprojective subrepresentations of $R_{n+1}$ and $L$ is a subrepresentation of $P_n\leq R_{n+1}$. After a look at the quotients $R_{n+1}/L$ and $P_n/L$ and using lemma~\ref{Lemma:Hom} one gets the following:
$$
\textrm{dim} Hom(L,R_{n+1}/L)^+=\textrm{dim}Hom(L,P_n/L)^++e_1
$$
and hence the cells of $Y'_0$ of dimension $k$ are in bijection with the cells of $Gr_{\mathbf{e}}(P_n)$ of dimension $k-e_1$ and \eqref{Eq:PY'e} holds.
\end{proof}

% Formulas\eqref{Eq:MainThmBetti}, \eqref{Eq:PoiPn} and \eqref{Eq:PoiIn} are also proved in an independent way by Szanto \cite{szanto-2009}.
\begin{corollary}
\begin{equation}
\chi(Gr_\mathbf{e}(R_n))={e_2\choose e_1}{n-e_1\choose e_2-e_1}.
\end{equation}
\begin{equation}
\chi(Gr_\mathbf{e}(P_n))={e_2-1\choose e_1}{n+1-e_1\choose e_2-e_1}
\end{equation}
\begin{equation}
\chi(Gr_\mathbf{e}(I_n))={e_2+1\choose e_1}{n-e_1\choose e_2-e_1}
\end{equation}
\end{corollary}
\begin{proof}
For a projective variety $X$, $\chi(X)=P_X(1)$ and $\chi(Gr_t(s))={s\choose t}$.
\end{proof}

\section{Application to cluster algebras}\label{Sec:ClusterAlgebras}
To a finite quiver $Q$ without loops and $2$--cycles  is associated a (coefficient--free) cluster algebra $\myAA_Q$ (\cite{FZI}, see also \cite{FZNotes} and \cite{keller-2008} for excellent surveys). This is a $\ZZ$--subalgebra of the field of rational functions in $n$ (= number of vertices of $Q$) variables generated by its \emph{cluster variables}. The \emph{cluster monomials} are monomials in cluster variables belonging to the same cluster.  A canonical basis of $\myAA_Q$ is a $\ZZ$--basis $\BB$ of it such that the positive linear combinations of elements of $\BB$ coincide with the semiring of positive elements of $\myAA_Q$, i.e. the elements which are positive Laurent polynomials in every cluster of $\myAA_Q$ (see \cite{Sherman}). The existence of such a basis has been proved only in a few cases: 
\begin{itemize}
 %\item if  $Q$ is of type $ADE$  then in \cite{hernandez}  and in \cite{Nakajima} prove that for a particular choice of ``coefficients'' there exists a monoidal categorification of $\myAA_Q$ and hence cluster monomials, i.e. monomials in cluster variables belonging to the same cluster, of $\myAA_Q$ form a basis whose elements are positive indecomposable, i.e. they cannot be written as sum of positive elements and hence they are a canonical basis.
\item when $Q$ is the Kronecker quiver, the canonical basis has been found by Sherman and Zelevinsky \cite{Sherman} (for every choice of the coefficients) and consists of cluster monomials together with extra elements $\{z_n:\;n\geq1\}$; 
\item when $Q$ is of type $A_2^{(1)}$, the canonical basis of $\myAA_Q$ has been found in \cite{CanonicalBasis} (for every choice of the coefficients) and consists of cluster monomials together with some extra elements $\{u_n:\geq 1\}$ possibly multiplied with some cluster variables.
\end{itemize}
Note that under the Caldero--Keller  bijection, cluster monomials correspond to rigid $Q$--representations and the quiver Grassmannians associated with rigid representations are smooth. In
the next two sections we give a geometric realization of the extra elements $\{z_n\}$, $\{u_n\}$.
 
\subsection{Type $A_1^{(1)}$}
Let $\mathcal{F}=\mathbf{Q}(x_1,x_2)$ be the field of rational functions in two independent variables $x_1$ and $x_2$ with rational coefficients. We define recursively elements $\{x_k|\; k\in\ZZ\}$ of $\myFF$ by:
$$
x_kx_{k+2}=x_{k+1}^2+1,\;\;\; k\in\ZZ.
$$
Let $\myAA$ be the $\ZZ$--subalgebra of $\myFF$ generated by the $x_k$ ($k\in\ZZ$). By \cite{Sherman} the algebra $\myAA$ is the coefficient--free cluster algebra associated with the Kronecker quiver. The couples $\{x_k,x_{k+1}\}$, $k\in\ZZ$, are free generating sets of $\myFF$ and they are called the \emph{clusters} of $\myAA$. Monomials $x_k^ax_{k+1}^b$, $a,b\geq0$, $k\in\ZZ$, are called the \emph{cluster monomials} of $\myAA$.

%By Fomin--Zelevinsky's Laurent phenomenon, every element of $\myAA$ can be written as a Laurent polynomial in every cluster of $\myAA$. The denominator of a cluster variable $x_k$, $k\neq 1,2$, in the cluster $\{x_1,x_2\}$ is a monomial of the form $x_1^{d_1(k)}x_2^{d_2(k)}$ for some non--negative integers $d_1(k), d_2(k)$.
%The vector $\mathbf{d}(k):=(d_1(k),d_2(k))$ is called the \emph{denominator vector} of $x_k$ in the cluster $\{x_1,x_2\}$. 
%Caldero and Keller \cite{CK1} have proved that the Caldero--Chapoton map (defined in the introduction) $M\mapsto CC(M)$ restricts to a bijection between indecomposable rigid $Q$--representations and cluster variables different from $x_1$ and $x_2$.

In \cite{CZ}, Caldero and Zelevinsky  have defined the rational function:
$$
s_n:=CC(R_n)%=\frac{\sum_{\mathbf{e}}\chi(Gr_\mathbf{e}(R_n))x_1^{n-e_2}x_2^{n}}{x_1^{n}x_2^{n}}
$$
for every $n\geq0$, where $R_n$ is a regular indecomposable $Q$--representation of dimension $(n,n)$. They have proved that the set
$$
\mathcal{S}:=\{\textrm{cluster monomials}\}\cup\{s_n:\;n\geq1\}
$$
is a $\ZZ$--basis of $\myAA$.

The canonical basis $\BB$ defined in the introduction and this basis are related by the following formula \cite{CZ}:
\begin{equation}\label{Eq:znsn}
z_n=s_n-s_{n-2}
\end{equation}
for $n\geq1$ and the convention that $s_{-k}=0$ for $k>0$.
\begin{theorem}
The element $z_n$ has the following Laurent expansion:
\begin{equation}\label{Eq:ZnLaurent}
z_n=\frac{\sum_{\mathbf{e}}\chi(Gr_\mathbf{e}(R_n)^{Sm})x_1^{2(n-e_2)}x_2^{2e_1}}{x_1^{n}x_2^{n}}
\end{equation}
where $Gr_\mathbf{e}(R_n)^{Sm}$ denotes the smooth part of $Gr_\mathbf{e}(R_n)$.
\end{theorem}
\begin{proof}
In view of theorem~\ref{Thm:Stratification}, $Gr_\mathbf{e}(R_n)^{Sm}=X_0\setminus X_1$ where $X_0=Gr_\mathbf{e}(R_n)$ and $X_1=Gr_{(e_1-1,e_2-1)}(R_{n-2})$ and hence:
$$
\chi(Gr_\mathbf{e}(R_n)^{Sm})=\chi(Gr_\mathbf{e}(R_n))-\chi(Gr_{(e_1-1,e_2-1)}(R_{n-2})).
$$
It is now easy to check that the right hand side of \eqref{Eq:ZnLaurent} satisfies \eqref{Eq:znsn}.
\end{proof}

\subsection{Type $A_2^{(1)}$}\label{Sec:A21}
\begin{figure}[hbt]
\begin{center}
\xymatrix{
         &&&       &         &                                &2\ar[rd]&\\
&&&&                         Q_2&1\ar[ur]\ar[rr]&             &3}
\end{center}
\caption{The quiver of type $A_2^{(1)}$}
\label{Fig:A21}
\end{figure}
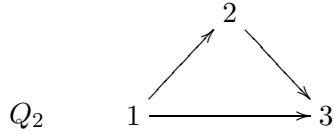

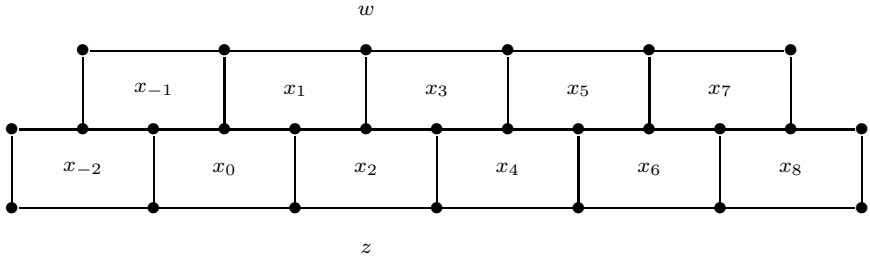
\begin{figure}
\xymatrix@C=20pt{
&&&&&&&&&&&&\\
&*{\bullet}\ar@{-}[d]\ar@{-}[rrrrrrrrrr]&&*{\bullet}\ar@{-}[d]&&*{\bullet}\ar@{-}[d]\ar@{}|{w}[u]&&*{\bullet}\ar@{-}[d]&&*{\bullet}\ar@{-}[d]&&*{\bullet}\ar@{-}[d]&\\
*{\bullet}\ar@{-}[d]\ar@{-}[rrrrrrrrrrrr]&*{\bullet}\ar@{}|{x_{-1}}[rru]&*{\bullet}\ar@{-}[d]&*{\bullet}\ar@{}|{x_1}[rru]&*{\bullet}\ar@{-}[d]&*{\bullet}\ar@{}|{x_3}[rru]&*{\bullet}\ar@{-}[d]&*{\bullet}\ar@{}|{x_5}[rru]&*{\bullet}\ar@{-}[d]&*{\bullet}\ar@{}|{x_7}[rru]&*{\bullet}\ar@{-}[d]&*{\bullet}&*{\bullet}\ar@{-}[d]\\
*{\bullet}\ar@{-}[rrrrrrrrrrrr]\ar@{}|{x_{-2}}[rru]&&*{\bullet}\ar@{}|{x_0}[rru]&&*{\bullet}\ar@{}|{x_2}[rru]\ar@{}|{z}[rrd]&&*{\bullet}\ar@{}|{x_4}[rru]&&*{\bullet}\ar@{}|{x_6}[rru]&&*{\bullet}\ar@{}|{x_8}[rru]&&*{\bullet}&\\
&&&&&&&&&&&&&
}
\caption{The exchange graph of a cluster algebra of type $A_2^{(1)}$}
\label{Fig:ExchA21}
\end{figure}

We now briefly recall the construction of the canonical basis of cluster algebras of type $A_2^{(1)}$ from \cite{CanonicalBasis}.  Let $\myFF =\QQ(x_1, x_2, x_3)$ be the field of rational functions in three (commuting) independent variables $x_1$, $x_2$ and $x_3$ with rational coefficients. Recursively define
elements $x_m\in \myFF$ for $m\in\ZZ$ by the relations
\begin{equation}
x_mx_{m+3} = x_{m+1}x_{m+2}+1.
\end{equation}
Define also the elements $w,z\in\myFF$ by
\begin{eqnarray}\label{Def:W}
w&=&\frac{x_1+x_3}{x_2}\\\label{Def:Z}
z&=&\frac{x_1x_2 +x_2x_3+1}{x_1x_3}
\end{eqnarray}
The (coefficient--free) cluster algebra $\myAA$ of type $A_2^{(1)}$ is the $\ZZ$--subalgebra of $\myFF$ generated by all the $x_m$, $w$ and $z$ (see also \cite[ Example 7.8]{FZI}). This is the cluster algebra associated with the affine quiver $Q_2$ of type $A_2^{(1)}$ shown in figure~\ref{Fig:A21}. 
The elements $x_m$, $m\in\ZZ$, $w$ and $z$ are the cluster variables of $\myAA$. The sets $\{x_m,x_{m+1},x_{m+2}\}$, $\{x_{2m}, z, x_{2m+2}\}$ and $\{x_{2m-1},w,x_{2m+1}\}$, $m\in\ZZ$, are the clusters of $\myAA$. The cluster monomials are monomials in cluster variables belonging to the same cluster. The exchange graph of $\myAA$ is the brick wall shown in figure~\ref{Fig:ExchA21}:  it has clusters as vertices and an edge between two vertices if the corresponding clusters share precisely two cluster variables. In this figure the cluster variables of a cluster label the regions surrounding its corresponding vertex.  
Define elements $u_n$, $n\geq0$, of $\myFF$ by the recursion:
\begin{equation}\label{Eq:UnDefSec}
\begin{array}{rcl}
u_0&=&2\\
u_1&=&zw-2\\
u_{n+1}&=&u_1 u_n-u_{n-1}\;\;\; n\geq1
\end{array}
\end{equation}
In \cite{CanonicalBasis} it is shown that the set 
\begin{equation}\label{Eq:CanBasA21}
\BB=\{\textrm{cluster monomials}\}\cup\{u_nw^k,\;u_nz^k:\;n\geq1,\;k\geq0\}
\end{equation}
is a $\ZZ$--basis of $\myAA$ such that the positive elements of $\myAA$ are precisely the non--negative linear combinations of elements of $\BB$.  We now realize the elements $u_n$ as image of the Caldero--Chapoton map. Recall that for a representation $M$ of the quiver $Q_2$ shown in figure~\ref{Fig:A21} the Caldero--Chapoton map $CC(M)$ is the following (see \cite{CZ}):
$$
CC(M)=\frac{\sum_{\mathbf{e}}\chi(Gr_\mathbf{e}(M))x_1^{d_2+d_3-e_2-e_3}x_2^{d_3-e_3+e_1}x_3^{e_1+e_2}}{x_1^{d_1}x_2^{d_2}x_3^{d_3}}
$$ 
where $(d_1,d_2,d_3)$ is the dimension vector of $M$. For every $n\geq 1$ let $R_{n,2}$ be the indecomposable regular $Q_2$--representation in an homogeneous tube, i.e.
$$
\xymatrix@R=8pt@C=10pt{
         &                                                     &k^n\ar[rd]^{=}&\\
R_{n,2}=&k^n\ar[ru]^{=}\ar[rr]_{J_n(0)}&                        &k^n}
$$ 
where $J_n(0)$ denotes the $n\times n$ indecomposable nilpotent Jordan block.
\begin{theorem}
For every $n\geq1$:
\begin{equation}\label{U=CC}
u_n=\frac{\sum_{\mathbf{e}}\chi(Gr_\mathbf{e}(R_{n,2})^{Sm})x_1^{2n-e_2-e_3}x_2^{n-e_3+e_1}x_3^{e_1+e_2}}{x_1^{n}x_2^{n}x_3^{n}}
\end{equation}
where $Gr_\mathbf{e}(R_{n,2})^{Sm}$ denotes the smooth part of the quiver Grassmannian $Gr_\mathbf{e}(R_{n,2})$.
\end{theorem}
\begin{proof}
Let $u'_n$ denote the right hand side of \eqref{U=CC}.
We consider the fibration of projective varieties
$$
\xymatrix@R=6pt{
f:Gr_{(e_1,e_2,e_3)}(R_{n,2})\ar@{->>}[r]&Gr_{(e_1,e_3)}(R_n)\\
(N_1,N_2,N_3)\ar@{|->}[r]&(N_1,N_3)
}
$$
whose fiber is $Gr_{(e_2-e_1)}(e_3-e_1)$. In particular
$$
Gr_{(e_1,e_2,e_3)}(R_{n,2})^{Sm}=f^{-1}\left(Gr_{(e_1,e_3)}(R_{n})^{Sm}\right).
$$
Using this it is straightforward to check that 
\begin{equation}
 u'_n=z_n(x_1w^{-1/2},x_3w^{-1/2})
\end{equation}
where the right hand side means to compute \eqref{Eq:ZnLaurent} with $x_1\mapsto x_1w^{-1/2}$ and $x_2\mapsto x_3w^{-1/2}$ and $w^{-1/2}$ is a formal variable whose square is $w^{-1}$. By definition \eqref{Eq:ZnDef} we have
$$
u'_1=z_1'(x_1w^{-1/2},x_3w^{-1/2})=\frac{x_1^2+x_3^2+w}{x_1x_2}=zw-2
$$
and $u'_{n+1}=u'_1u'_n-u'_{n-1}$ and hence $u'_n=u_n$ for every $n\geq1$.
\end{proof}

\section{Conclusions}
The study of quiver Grassmannians from a geometric point of view has given interesting information about cluster algebras. After this study it is natural to define  the following slight modification of the Caldero--Chapoton map: for $k\geq0$
\begin{equation}\label{Eq:CCModification}
CC^{(k)}(M):=\frac{\sum_{\mathbf{e}}\chi(Gr_\mathbf{e}(M)^{(k)})x_1^{2(d_2-e_2)}x_2^{2e_1}}{x_1^{d_1}x_2^{d_2}}
\end{equation}
where
$$
Gr_\mathbf{e}(M)^{(k)}:=\{N\in Gr_\mathbf{e}(M)|\; \textrm{dim}Ext^1(N,M/N)=k\}.
$$
So $CC^{(0)}(M)$ is defined by the smooth part of the quiver Grassmannians associated with $M$. In particular for a rigid representation $M$ we have $CC^{(0)}(M)=CC(M)$. 
% and in view of theorem~\ref{Eq:ZnLaurent} the Sherman--Zelevinsky canonical basis $\BB$ consists section~\ref{Sec:ClusterAlgebras} looks like:
% $$
% \BB=\{CC^{(0)}(M)|\; M\textrm{ rigid }\}\cup\{CC^{(0)}(R_n)|\; n\geq1\}\cup\{x_1^ax_2^b|\;a,b\geq0\}.
% $$
Moreover it is easy to check that if $Ext^1(M,N)=Ext^1(N,M)=0$ then 
$CC^0(M\oplus N)=CC^0(M)CC^0(N)$. In particular the elements $u_nw$ and $u_nz$ of the basis \eqref{Eq:CanBasA21} are $CC^0(R_{n,2}\oplus W)$ and $CC^0(R_{n,2}\oplus Z)$ where $W$ and $Z$ are the indecomposable regular rigid $Q_2$--representations.

It would be interesting to study the relation between the nice construction due to Dupont \cite{dupont-2009} of transverse quiver Grassmannians with the maps $CC^{(k)}$. 
%We find that the transverse quiver Grassmannians of indecomposable regular homogeneous representations coincide with the smooth part of the quiver Grassmannians both in type $A_1^{(1)}$ and $A_2^{(1)}$. But in tubes of higher rank the two notions slightly differ.

We hope that the technique developed here can help in the determination of a ``canonical basis'' of cluster algebras of higher rank which is still an open problem. 
%The only solved cases are type $A_1^{(1)}$ \cite{Sherman} and $A_2^{(1)}$ \cite{CanonicalBasis}.

\subsection*{Acknowledgements}
The first author thanks the organizers of the ``Advanced School and Conference on Homological and Geometrical Methods in Representation Theory" which took place in Trieste (Italy) at the beginning of 2010 for the opportunity to present this project in front of several experts and to spend a fruitful time there.

We are grateful to Andrei Zelevinsky for having pointed out the formulas found by Szanto in \cite{szanto-2009}.

We thank Gr\`egoire Dupont for useful comments about transverse quiver Grassmannians.

\bibliographystyle{plain}
\bibliography{BibliografiaGiovanni}

\begin{thebibliography}{10}

\bibitem{ASS}
Ibrahim Assem, Daniel Simson, and Andrzej Skowro{\'n}ski.
\newblock {\em Elements of the representation theory of associative algebras.
  {V}ol. 1}, volume~65 of {\em London Mathematical Society Student Texts}.
\newblock Cambridge University Press, Cambridge, 2006.
\newblock Techniques of representation theory.

\bibitem{BFZIII}
Arkady Berenstein, Sergey Fomin, and Andrei Zelevinsky.
\newblock Cluster algebras. {III}. {U}pper bounds and double {B}ruhat cells.
\newblock {\em Duke Math. J.}, 126(1):1--52, 2005.

\bibitem{BB}
A.~Bia{\l}ynicki-Birula.
\newblock Some theorems on actions of algebraic groups.
\newblock {\em Ann. of Math. (2)}, 98:480--497, 1973.

\bibitem{CC}
Philippe Caldero and Fr{\'e}d{\'e}ric Chapoton.
\newblock Cluster algebras as {H}all algebras of quiver representations.
\newblock {\em Comment. Math. Helv.}, 81(3):595--616, 2006.

\bibitem{CK2}
Philippe Caldero and Bernhard Keller.
\newblock From triangulated categories to cluster algebras. {II}.
\newblock {\em Ann. Sci. \'Ecole Norm. Sup. (4)}, 39(6):983--1009, 2006.

\bibitem{CK1}
Philippe Caldero and Bernhard Keller.
\newblock From triangulated categories to cluster algebras.
\newblock {\em Invent. Math.}, 172(1):169--211, 2008.

\bibitem{CR}
Philippe Caldero and Markus Reineke.
\newblock On the quiver {G}rassmannian in the acyclic case.
\newblock {\em J. Pure Appl. Algebra}, 212(11):2369--2380, 2008.

\bibitem{CZ}
Philippe Caldero and Andrei Zelevinsky.
\newblock Laurent expansions in cluster algebras via quiver representations.
\newblock {\em Mosc. Math. J.}, 6(3):411--429, 2006.

\bibitem{Chriss}
Neil Chriss and Victor Ginzburg.
\newblock {\em Representation theory and complex geometry}.
\newblock Birkh\"auser Boston Inc., Boston, MA, 1997.

\bibitem{CrawleyTree}
W.~W. Crawley-Boevey.
\newblock Maps between representations of zero-relation algebras.
\newblock {\em J. Algebra}, 126(2):259--263, 1989.

\bibitem{DWZII}
Harm Derksen, Jerzy Weyman, and Andrei Zelevinsky.
\newblock Quivers with potentials and their representations ii: Applications to
  cluster algebras, 2009.

\bibitem{dupont-2009}
G.~Dupont.
\newblock Transverse quiver grassmannians and bases in affine cluster algebras,
  2009.
\newblock arXiv.org:0910.5494.

\bibitem{FZI}
Sergey Fomin and Andrei Zelevinsky.
\newblock Cluster algebras. {I}. {F}oundations.
\newblock {\em J. Amer. Math. Soc.}, 15(2):497--529 (electronic), 2002.

\bibitem{FZII}
Sergey Fomin and Andrei Zelevinsky.
\newblock Cluster algebras. {II}. {F}inite type classification.
\newblock {\em Invent. Math.}, 154(1):63--121, 2003.

\bibitem{FZNotes}
Sergey Fomin and Andrei Zelevinsky.
\newblock Cluster algebras: notes for the {CDM}-03 conference.
\newblock In {\em Current developments in mathematics, 2003}, pages 1--34. Int.
  Press, Somerville, MA, 2003.

\bibitem{FZIV}
Sergey Fomin and Andrei Zelevinsky.
\newblock Cluster algebras. {IV}. {C}oefficients.
\newblock {\em Compos. Math.}, 143(1):112--164, 2007.

\bibitem{hernandez}
David Hernandez and Bernard Leclerc.
\newblock Cluster algebras and quantum affine algebras.
\newblock arXiv.org:0903.1452.

\bibitem{CanonicalBasis}
Giovanni~Cerulli Irelli.
\newblock Canonically positive basis of cluster algebras of type $a_2^{(1)}$,
  2009.
\newblock arXiv.org:0904.2543.

\bibitem{String}
Giovanni~Cerulli Irelli.
\newblock Quiver grassmannians associated with string modules, 2009.
\newblock arXiv.org:0910.2592.

\bibitem{keller-2008}
Bernhard Keller.
\newblock Cluster algebras, quiver representations and triangulated categories,
  2008.
\newblock arXiv.org:0807.1960.

\bibitem{Kronecker}
Leopold Kronecker.
\newblock Algebraische reduction der schaaren bilinearer formen.
\newblock {\em Sitzungsberichte Akad. Berlin}, pages 1225--1237, 1890.

\bibitem{Nakajima}
Hiraku Nakajima.
\newblock Quiver varieties and cluster algebras, 2009.
\newblock arXiv.org:0905.0002.

\bibitem{RingelExceptional}
Claus~Michael Ringel.
\newblock Exceptional modules are tree modules.
\newblock {\em Proceedings of the {S}ixth {C}onference of the {I}nternational
  {L}inear {A}lgebra {S}ociety ({C}hemnitz}.

\bibitem{Sherman}
Paul Sherman and Andrei Zelevinsky.
\newblock Positivity and canonical bases in rank 2 cluster algebras of finite
  and affine types.
\newblock {\em Mosc. Math. J.}, 4(4):947--974, 982, 2004.

\bibitem{szanto-2009}
Csaba Szanto.
\newblock On the cardinalities of kronecker quiver grassmannians, 2009.
\newblock arXiv.org:0903.1928.

\end{thebibliography}

\end{document}